\documentclass[12pt]{article}

\usepackage{hyperref}

\usepackage{amsmath}
\usepackage{amssymb}
\usepackage{amsfonts}
\usepackage{amsopn}
\usepackage[all]{xy}
\usepackage{amsthm}
\usepackage{amscd}

\swapnumbers
\theoremstyle{plain}
\newtheorem{thm}{Theorem}[section]
\newtheorem{lem}[thm]{Lemma}

\newtheorem{prop}[thm]{Proposition}

\theoremstyle{definition}
\newtheorem{ntt}[thm]{}
\newtheorem{rem}[thm]{Remark}

\newcommand{\af}{\mathbb{A}}   
\newcommand{\pl}{\mathbb{P}}   
\newcommand{\zz}{\mathbb{Z}}   
\newcommand{\hh}{\mathbb{H}}   

\newcommand{\mot}{\mathcal{M}_A}             
\newcommand{\smp}{\mathcal{S}m\mathcal{P}}  
\newcommand{\ab}{\mathcal{A}b}               
\newcommand{\corr}{\mathcal{C}or_A}          
\newcommand{\smpr}{\mathcal{S}m\mathcal{P}roj}            

\newcommand{\ra}{\rightarrow}  
\newcommand{\hra}{\hookrightarrow} 
\newcommand{\xra}[1]{\xrightarrow{#1}} 

\DeclareMathOperator{\spec}{Spec}   
\DeclareMathOperator{\codim}{codim} 

\newcommand{\id}{\mathrm{id}}       

\title{Oriented Cohomology and Motivic Decompositions of Relative Cellular Spaces}
\author{A.~Nenashev, K.~Zainoulline}
\date{September 24, 2004}

\begin{document}

\maketitle

\begin{abstract}
For an oriented cohomology theory $A$ and a relative cellular space $X$,
we decompose the $A$-motive of $X$ into a direct sum of twisted motives 
of the base spaces. We also obtain respective decompositions of the 
$A$-cohomology of $X$.
Applying them, one can compute
$A(X)$, where $X$ is an isotropic projective homogeneous variety 
and
$A$ means algebraic $K$-theory, motivic cohomology or
algebraic cobordism $MGL$.

\ 


\end{abstract}

\section{Introduction}
A cellular decomposition of a smooth projective variety $X$ is a filtration by 
closed subvarieties
$X=X_0\supset X_1\supset\ldots\supset X_n\supset X_{n+1}=\emptyset$
along with an affine bundle structure on each complement $X_i\smallsetminus X_{i+1}$ over a 
certain smooth projective base $Y_i$. 
Our first goal is to compute $A(X)$ in terms of $A(Y_i)$, 
where $A$ 
is an arbitrary oriented cohomology theory in the sense of \cite{Pa}. 
Namely, we prove that
\begin{equation}\label{cdecomp}
A^*(X)\cong \bigoplus_{i=0}^n A^{*-c_i}(Y_i),
\end{equation}
where $c_i$ is the codimension of $X_i$ in $X$ (Thm.~\ref{decthm}). 
Note that the shift of degree by $-c_i$ in (\ref{cdecomp}) holds under 
the assumption that $A$ is $\zz$-graded and 
push-forwards in $A$ raise the degree by the 
codimension of the morphism. 
If it is twice the codimension, e.g. for the singular cohomology
of complex varieties, then we get $-2c_i$
in (\ref{cdecomp}). 
If $A$ is bigraded, which is the case for many theories,  
then the formula must also be changed accordingly, 
e.g. for the algebraic cobordism we have
$$
MGL^{2*,*}(X)\cong \bigoplus_{i=0}^nMGL^{2*-2c_i,*-c_i}(Y_i).
$$

Our second goal is to introduce the category of (graded) 
$A$-correspondences $\corr$. The cohomology functor $A$ factors
through correspondences, and we show that the isomorphism (\ref{cdecomp})
comes in fact from a certain isomorphism in the category $\corr$.
Taking the transposed isomorphism in $\corr$ and translating it 
back in terms of $A$, we get a dual formula for $A(X)$:
\begin{equation}\label{rdecomp}
A^*(X)\cong \bigoplus_{i=0}^n A^{*-r_i}(Y_i),
\end{equation}
where $r_i$ denotes the rank of $(X_i\smallsetminus X_{i+1})\ra Y_i$
as an affine bundle (Thm.~\ref{corrdcp}). 

Finally, following Grothendieck \cite{Gr} and Manin \cite{Ma} 
we define the category
of $A$-motives $\mot$ by taking the pseudo-abelian completion of $\corr$.
Then formula (\ref{rdecomp}) can be interpreted in $\mot$ as a motivic 
decomposition of the cellular variety $X$:
\begin{equation}\label{motdecomp}
X\cong\bigoplus_{i=0}^n Y_i\otimes L^{\otimes r_i},
\end{equation}
where $L$ is the Lefschetz motive (Thm. \ref{mainthm}).

Our research was inspired by the work of N.~Karpenko
\cite{Ka} who calculated the Chow groups of relative cellular varieties.
We wanted to extend Karpenko's results, using his techniques wherever
possible, to arbitrary oriented cohomology theories, having in mind
the algebraic cobordism $MGL$ as a most important example. 
Note that Chow groups provide a special example of an oriented theory 
in that they have some nice extra properties that one cannot expect
of other examples. In particular, there are tools enabling one to work
with the Chow groups of non-regular varieties, which is used essentially 
in Karpenko's approach, where some of the varieties arising in the arguments
do not need to be smooth, and push-forwards along morphisms of such
varieties are considered. This makes a straightforward extension of the
methods of \cite{Ka} to the general case impossible since the
machinery of oriented cohomology theories, producing push-forwards \cite{Pa1},
has been developed for smooth varieties only.

For this reason, the initial version of the paper was written for cellular
filtrations in which all terms $X_i$ were smooth varieties, but even in
this case essential modifications of the arguments were already necessary.
Briefly speaking, in order to construct a splitting of a certain
localization sequence used in the proof of (\ref{cdecomp}), we have to refer
to a result of Hironaka on the elimination of points of indeterminacy \cite{Hi}.
It is there that we restrict
ourselves to ground fields of characteristic $0$.
Having this splitting and assuming that all terms $X_i$
are smooth varieties, we immediately get the decompositions (\ref{cdecomp}),
(\ref{rdecomp}) and (\ref{motdecomp}).   


Later A. Merkurjev pointed out to us that the assumption of smoothness of
the $X_i$ could be removed if one works with cohomology groups with support
rather than replacing them with absolute cohomology by means of Gysin operators 
(recall that supports are not assumed to be smooth in \cite{Pa1}). 
To implement this idea, one must first develop a theory of push-forwards
for cohomology with support (as opposed to the absolute case), which has not 
been done in \cite{Pa1} or elsewhere. We show how to define such
push-forwards. However, we do not try
to make a complete theory out of our definition,
rather we restrict ourselves to what we really need.


The paper is organized as follows. First,
we review the notion of an oriented cohomology theory (Sect.~\ref{orcohth})
and then define push-forwards for cohomology with support, proving 
only the expected properties that are necessary for our purposes 
(Thm.~\ref{relgysin}).
Then (Sect.~\ref{celldeccoh}) we obtain decomposition (\ref{cdecomp}). 
Next (Sect.~\ref{cellcorrdec}), 
we introduce the category of $A$-correspondences
and prove decomposition (\ref{rdecomp}) .
In section \ref{amot} we define the category of $A$-motives
and prove motivic decomposition (\ref{motdecomp}).
Finally, we provide some examples of relative cellular spaces
with smooth filtrations and compute their cohomology.

\paragraph{Conventions, notation, terminology.}
We fix a field $k$ and consider algebraic varieties over it assuming
that they are quasi-projective. Whenever we speak of the {\it codimension
of a morphism}, it means that the morphism is assumed to be 
equicodimensional. However all the results can be extended to the case
of non-equicodimensional morphisms by writing separate formulas on
components. 

We are making intensive use of the results of I. Panin presented in 
\cite{Pa} and \cite{Pa1}, where no degree indexation for cohomology 
groups is being displayed and the cohomology theories take their values 
in the category $\ab$ of abelian groups, rather than in $\ab^{\zz}$ or 
$\ab^{\zz/2}$. For this reason we sometimes `forget' to write upper 
indices in cohomology, though our theories are $\zz$-graded; 
we hope this will not lead to a confusion.

\paragraph{Acknowledgements} 
Both authors thank the University of Bielefeld
for hospitality and Alexander von Humboldt Stiftung for support.
The last author also thanks the RTN-Network HPRN-CT-2002-00287
for support. We are grateful to Manuel Ojanguren who has kindly 
invited us to spend a week at the EPFL, where 
we had an opportunity to discuss our work with him and I. Panin.
We are indebted to B.~Calm\`es, V.~Guletskii, A.~Merkurjev, I.~Panin, 
A.~Pukhlikov, M.~Rost and A.~Vishik for stimulating discussions.

\section{Oriented cohomology theories}\label{orcohth}
In the present section we recall the notion of an oriented cohomology
theory $A$ as defined in \cite{Pa} and \cite{Pa1}.

\begin{ntt}
Let $A^*:\smp^{op} \ra \ab^{\zz}$ be a contravariant functor 
from the {\it category of smooth pairs}
over a field $k$ to the category of $\zz$-graded abelian groups,
and let $A^i$ denote its $i$-th graded component.
Objects of $\smp$ are pairs $(X,U)$ consisting of a smooth quasi-projective
variety $X$ over $k$ and an open subset $U\subset X$.
A morphism $f$ from $(X,U)$ to $(X',U')$ 
is a morphism of pairs, i.e., an $f:X\ra X'$ such that $f(U)\subset U'$. 
The category of smooth quasi-projective varieties over $k$ embeds into $\smp$
by $X\mapsto (X,\emptyset)$; we will denote this object of $\smp$ 
simply by $X$.
By $pt=\spec k$ we denote the final object of $\smp$.

Let $Z$ be a closed subset of $X$.
We write $A_Z^i(X)$ for $A^i(X,X\smallsetminus Z)$ having in mind the notation
used for cohomology with support.  
We also write $A^i(X)$ for $A^i(X,\emptyset)$.
For a given morphism $f:(X,X\smallsetminus Z)\ra (X',X'\smallsetminus Z')$ the induced map 
$f^*:A^i_{Z'}(X') \ra A^i_Z(X)$ is called a {\it pull-back}.
Note that pull-backs preserve the grading.
\end{ntt}

\begin{ntt}
{\it A ring cohomology theory} is a contravariant functor 
$A^*:\smp^{op} \ra\ab^{\zz}$ together with a functorial 
morphism of degree $+1$
$$
\partial: A^*(U) \ra A^{*+1}(X,U)$$ 
and a cup-product
$$
\cup: A_Z^i(X)\times A_{Z'}^j(X) \ra A_{Z\cap Z'}^{i+j}(X)
$$ 
with usual properties
(see \cite[Sect.~2]{Pa}). 
Observe that the cup-product turns $A^*(X)$ into
a $\zz$-graded commutative ring, i.e., 
$\alpha\cup\beta=(-1)^{\deg(\alpha)\deg(\beta)}\beta\cup\alpha$, 
for homogeneous $\alpha$, $\beta \in A^*(X)$.
\end{ntt}

In the sequel we will refer to the following properties of a
ring cohomology theory $A^*$.
\begin{ntt}[Localization sequence]\label{locseq}
Let $Y$ be a closed subset of a smooth variety $X$ and $Z$ be a closed
subset of $Y$. Then there is an exact sequence of $A^*(pt)$-modules
(see \cite[2.2.3]{Pa})
$$
\ldots\ra A_Z^*(X)\ra A_Y^*(X) \ra A_{Y-Z}^*(X-Z) 
\xra{\partial} A_Z^{*+1}(X) \ra \ldots
$$
\end{ntt}

\begin{ntt}[Strong Homotopy Invariance]\label{hominvar}
Let $p: E\ra X$ be an affine bundle over a smooth variety $X$
and $Z\hra X$ be a closed subset.
Then the pull-back $p^*:A^*_Z(X)\ra A^*_{p^{-1}(Z)}(E)$ is an isomorphism
(see \cite[2.2.6]{Pa}).
\end{ntt}

\begin{ntt}\label{orcohprop}
{\it An oriented cohomology theory} 
is a ring cohomology theory $A^*$
together with a rule that assigns to each vector bundle $E/X$ of rank $r$
and a closed subset $Z\hra X$ 
an operator
$th^E_Z: A_Z^*(X) \ra A_Z^{*+r}(E)$
which is a two-sided $A^*(X)$-module isomorphism 
and satisfies invariance, base change and additivity properties
(see \cite[Def.~3.1]{Pa}). This operator will be referred to as a 
{\it Thom isomorphism}.

An equivalent way to make a theory oriented is to endow it with 
first Chern classes for line bundles, i.e., to assign to each
line bundle $L$ over a smooth $X$ an element 
$c(L)\in A^1(X)$ with suitable properties (see \cite{Pa}). 
\end{ntt}

\begin{ntt}
Among examples of oriented cohomology theories are 
higher Chow groups defined by Bloch, motivic cohomology, 
$K$-Theory, semi-topological $K$-theory introduced by Friedlander,
\'etale cohomology with coefficients in $\mu_n$, and the
algebraic cobordism $MGL$
(see \cite[Sect.~3.8]{Pa} for more examples).

Observe that the notion of an oriented 
cohomology theory introduced in \cite{LM}
is more general than that of \cite{Pa}. For there are 
no localization sequences \ref{locseq} on the list of axioms in \cite{LM}.
Hence, we cannot say that the algebraic cobordism theory 
$\Omega$ defined in \cite{LM}
is an example of an oriented theory in the sense of \ref{orcohprop}.
\end{ntt}

An oriented cohomology theory satisfies the following important properties
which will be used in the sequel.

\begin{ntt}[Projective Bundle Theorem]\label{projbun} 
Let $E/X$ be a vector
bundle of rank $n+1$ over a smooth variety $X$ and $Z\hra X$ a closed subset.
Then the map
$$
(1,\xi,\ldots,\xi^n)\cup -: 
\bigoplus_{i=0}^n A^{*-i}_Z(X) \ra A_{\pl(E_Z)}^*(\pl(E)),
$$
is an isomorphism, where $\xi=c(\mathcal{O}_E(-1))$ and $E_Z = E|_Z$ 
is the restriction of $E$ to $Z$ (see \cite[Cor.~3.17]{Pa}).
\end{ntt}

\begin{ntt}[Integration structure]\label{integr}
According to \cite[Thm.~4.1.4]{Pa1}, an oriented cohomology theory $A^*$
has a unique integration structure. This means that
for any projective morphism of smooth varieties $f:Y\ra X$ of codimension $c$
there is given a two-sided $A^*(X)$-module operator
$$
f_*:A^*(Y)\ra A^{*+c}(X)
$$ 
called a {\it push-forward} which satisfies certain properties 
(see \cite[Def.~4.1.2]{Pa1}).
\end{ntt}

\begin{ntt}[Projection Formula]\label{projform}
Let $f:Y\ra X$ be a projective morphism of smooth varieties.
Then for any  $\alpha\in A(Y)$ and $\beta\in A(X)$
$$
f_*(\alpha\cup f^*(\beta))=f_*(\alpha)\cup \beta
.
$$
\end{ntt}

\begin{ntt}[Push-forwards for closed embeddings] \label{gysin}
For a closed embedding of smooth varieties $\imath:Y\ra X$ of codimension $c$ and 
a closed subset $Z\hra Y$, not necessarily smooth, the $A^*(X)$-module operator
$\imath_{th}$ is defined in \cite{Pa} as the composition 
$$
\imath_{th}: A_Z^*(Y) \xra{th^N_Z} A_Z^{*+c}(N) \xra{d_Z(X,Y)}  A_Z^{*+c}(X)\, .
$$
Here $th^N_Z$ is the Thom isomorphism of \ref{orcohprop}, $N=N_Y X$ is the 
normal bundle to $Y$ in $X$, and $d_Z(X,Y)$ is the deformation to the
normal cone isomorphism for the pair $Y\hra X$ with support in $Z$
(see \cite[Thm.~2.2]{Pa1} for more details). Thus $\imath_{th}$ is an isomorphism as well.

Observe that the composition
$$
A^*(Y)\xra{\imath_{th}} A_Y^{*+c}(X) \ra A^{*+c}(X),
$$
where the last map is an extension of support, is the push-forward
$\imath_*$ (Gysin map) induced by the embedding $\imath$.
\end{ntt}

\section{Push-forwards with support}

\begin{ntt}
As it was mentioned in the introduction, the theory of push-forwards 
for cohomology with support has not been so far developed.
Partial results were obtained in \cite{Sm}, where the case of closed
embeddings was treated, under the name of {\it local push-forwards}.
It can be expected to be a long exercise to extend the methods of
\cite{Pa1}, where push-forwards along arbitrary projective morphisms
are constructed for absolute cohomology groups, in order to get
a theory of push-forwards for cohomology groups with support. 
Avoiding any comments as to how long and technically complicated such 
an exercise could be, we want to point out that it is clear how to 
define such push-forwards.

Let $f: X\to X'$ be a projective morphism of codimension $d$
of smooth varieties over $k$. Let $Z\hra X$ and $Z'\hra X'$ be 
closed subsets such that $f(Z)\hra Z'$.
Let $f=p\circ \imath$ be a decomposition of $f$ into a closed embedding
 $\imath: X\hra X'\times\pl^m$
followed by the projection $p:X'\times\pl^m \to X'$.
Define a $A^*(X)$-module homomorphism
$f_*: A^*_Z(X) \ra A^{*+d}_{Z'}(X')$ as the composition
\begin{equation}\label{pfsup}
A_Z(X) \xra{\imath_{th}} A_{\imath(Z)}(X'\times\pl^m) 
\to A_{Z'\times \pl^m}(X'\times\pl^m)
\xra{p_*} A_{Z'}(X'),
\end{equation}
where $\imath_{th}$ is the isomorphism of \ref{gysin},
the middle map is an extension of support, 
and the last map $p_*$ is defined in the proof of Theorem \ref{relgysin}.

In order to show that the assignment $f\mapsto f_*$ indeed provides
a theory of push-forwards (also known as {\it integration} in the 
terminology of \cite{Pa1})
for cohomology with support,
one must prove that these $f_*$-s
do not depend on the choice of decompositions $f=p\circ \imath$
and that they satisfy some standard properties of push-forwards 
(see \cite[Sect.4]{Pa1}). 
The proofs of these properties are still lacking in full generality.
However, in the present paper we will only use push-forwards with support
in the situation described by the following theorem.
\end{ntt}

\begin{thm}\label{relgysin}
There exists a rule which assigns to 
\begin{itemize}
\item
every equicodimensional projective morphism $f: W\to X$ of smooth varieties
endowed with a decomposition into a closed embedding $\imath: W\to X\times \pl^m$
followed by the projection $p: X\times\pl^m \to X$, and
\item
a closed subset $Z$ in $X$ satisfying $f(W)\subset Z$
\end{itemize}
a map $\overline{f_*}: A^*(W) \ra A^{*+c}_{Z}(X)$, 
where $c$ is the codimension of $f$. These maps have the 
following properties
\begin{itemize}
\item[\rm (i)]
For any open embedding $u: U\hra X$ with $U\cap Z\neq\emptyset$,
the diagram
\begin{equation}\label{diagrpush}
\xymatrix{
A_Z^{*+c}(X) \ar[r]^{u^*}& A_{Z_U}^{*+c}(U)\\
A^*(W) \ar[u]^{\bar{f}_*}  \ar[r]^{u^*} & A^*(W_U) \ar[u]_{\overline{({f}_U)}_*}
}
\end{equation}
commutes, where $Z_U=Z\times_X U$, $W_U=W\times_X U$, and we assume that 
the decomposition for $f_U$ is induced by the decomposition for $f$.
\item[\rm (ii)] 
The composition $A(W) \xra{\overline{f_*}} A_Z(X) \ra A(X)$ is the usual push-forward $f_*$
of \cite{Pa1} (here the last map is an extension of support).
Thus it does not depend on the choice of a decomposition $f=p\circ\imath$.
\item[\rm (iii)] 
If $f$ is a closed embedding and $Z=f(W)$, then $\overline{f_*}=f_{th}$ (see 2.10) for any
choice of the decomposition.
\end{itemize}
\end{thm}

\begin{proof}
This is a particular case of (\ref{pfsup}).
Define $\overline{f_*}$ to be the composition
$$
A(W) \xra{\imath_{th}} A_{W}(X\times\pl^m) \ra A_{Z\times \pl^m}(X\times\pl^m)
\xra{p_*} A_Z(X),
$$
where $\imath_{th}$ is the isomorphism of \ref{gysin},
the middle map is an extension of support, 
and the last map $p_*$ is defined as follows.

According to \ref{projbun}, an element $a\in A_{Z\times\pl^m}(X\times \pl^m)$
can be uniquely written as a sum 
$$a=a_0+\xi\cup a_1+\ldots  \xi^m\cup a_m,$$
where $a_i\in A_Z(X)$ and $\xi=c(\mathcal{O}_E(-1))$.
Define a $A_Z(X)$-linear map $p_*$ by setting 
$$p_*(\xi^i)=g^*([\pl^{m-i}]_\omega)\in A(X),$$ where $g^*$ is the pull-back induced
by the structural morphism $g:X\ra pt$,
$[\pl^{m-i}]_\omega$ is the image of the class of the projective space
by means of the induced map $\mathbb{L}\ra A(pt)$ from the Lazard ring
to the ring of coefficients of the cohomology theory $A$ (see \cite[4.3.1]{Pa1}).
In other words, the definition of $p_*$ `with support' is the same as
in the absolute case.

Note that the isomorphisms $\imath_{th}$ 
respect base changes via open embeddings
according to \cite[4.4.4]{Pa1}, and so do  the push-forwards $p_*$.
The latter is proved in \cite[4.5]{Pa1} in the absolute case,
but the proof equally works `with support', which proves 
assertion (i) of the theorem.

The diagram
$$
\xymatrix{
A_{Z\times\pl^m}(X\times\pl^m)\ar[r]^-{p_*} \ar[d] & A_Z(X)\ar[d] \\
A(X\times \pl^m) \ar[r]^-{p_*} & A(X)
}
$$
commutes by the very definition of both $p_*$-s. It follows that
the composition $A(W) \xra{\overline{f_*}} A_Z(X) \ra A(X)$ coincides with
$A(W) \xra{\imath_*} A(X\times \pl^m) \xra{p_*} A(X)$, which
is the definition of $f_*$ in \cite{Pa1}. This proves assertion (ii). 

Let $f: W\to X$ be a closed embedding such that $Z=f(W)$.
It remains to check that $\overline{f_*}$ does not depend on the choice
of a decomposition $f=p\circ\imath$.

Consider the Cartesian square
$$
\xymatrix{
W\times \pl^m \ar[r]^{f\times\id}\ar[d]^{q}& X\times \pl^m\ar[d]^{p} \\
W \ar[r]^{f}\ar[ur]^{\imath} \ar@/^/[u]^s& X
}
$$
(here $q$ denotes the projection). 
Clearly $\imath$ can be written as $\imath = (f,g)$, with $g: W \ra\pl^m$.
Then the closed embedding $s = (\id,g): W \ra W\times \pl^m$ is a section of
$q$ satisfying  $\imath =(f\times\id)\circ s$.

Consider the diagram
$$
\xymatrix{
A(W)\ar[r]^-{s_{th}}\ar[rd]^-{s_*}\ar[rdd]_-{\id}  & A_{s(W)}(W\times\pl^m)\ar[r]^-{(f\times\id)_{th}} \ar[d]_-{e} & 
A_{\imath(W)}(X\times \pl^m) \ar[d]^-e \\
      & A(W\times\pl^m)\ar[r]^-{(f\times\id)_{th}} \ar[d]_-{q_*} & 
A_{f(W)\times\pl^m}(X\times\pl^m) \ar[d]^-{p_*} \\
      & A(W) \ar[r]^-{f_{th}}           & A_{f(W)}(X)
}
$$
where $e$ is an extension of support, $s_*=e\circ s_{th}$ is the usual
push-forward for the closed embedding $s$, and $q_*\circ s_*=\id$ by 
\cite[4.6.(56)]{Pa1}. Commutativity of the upper square 
amounts to the following two assertions.
\begin{itemize}
\item
Thom isomorphisms commute with the extension of support maps $e$ (in fact, 
they are compatible with any pull-backs, see \cite{Pa1}).
\item
The deformation to the normal cone isomorphisms commute with the $e$-s,
which is stated in \cite{Pa}; see also \cite[3.2.2]{Sm} for a more accurate account
in the situation with support. 
\end{itemize}
The bottom square commutes by \cite[4.5]{Pa1} (see the proof
of property 4 of Quillen operators which equally works `with support').
The outer contour of the diagram can be now simplified to
$$
\xymatrix{
A(W) \ar[r]^-{\imath_{th}} \ar[rd]_-{f_{th}}& A_{\imath(W)}(X\times\pl^m)\ar[d]^-{p_*\circ e} \\
     & A_{f(W)}(X)
}
$$
(here we use the property that the maps $(-)_{th}$ are compatible with compositions, 
see \cite[3.2.2]{Sm} or \cite{Ne1}), and we are done.
\end{proof}

\section{Oriented cohomology of a relative cellular space}\label{celldeccoh}

\begin{ntt}
Recall that a scheme $E$ together with a morphism $p: E\ra X$
is an {\it affine bundle} of rank $r$ over $X$ if $X$ admits an open covering
by $U_j$ and isomorphisms 
$p^{-1}(U_j)\cong U_j\times \af^r$ over $U_j$
such that the transition isomorphisms $(U_i\cap U_j)\times \af^r\ra (U_i\cap U_j)\times \af^r$ are affine transformations in the obvious sense.
\end{ntt}

\begin{ntt}
Let $X$ be a smooth projective variety.
We call $X$ a {\it relative cellular space} if there exists a 
finite decreasing filtration by closed subsets 
$$
X=X_0\supset X_{1}\supset \ldots \supset X_{n}\supset \emptyset
$$
such that each complement $E_i=X_i\smallsetminus X_{i+1}$ 
is an affine bundle of a constant rank $r_i$ over a smooth
projective variety $Y_i$ (we assume $Y_n=X_n$). 
We denote by $p_i:E_i\ra Y_i$ the projection morphism.
\end{ntt}

\begin{ntt} Among examples of relative cellular spaces are
isotropic projective homogeneous varieties $H=G/P$, 
where $G$ is a linear algebraic group
(not necessarily split) and $P$ a parabolic subgroup (see \cite{CGM}), and,
more generally, smooth projective varieties
equipped with an action of the multiplicative
group $\mathbb{G}_m$ (see \cite[Thm.~3.4]{Br}). 
\end{ntt}

\begin{thm}\label{decthm}
Let $A^*$ be an oriented cohomology theory over a field of
characteristic 0. Let $X$ be a relative cellular
space with a filtration 
$X=X_0\supset X_1 \supset \ldots \supset X_n$. Assume that
all $X_i$ are equi-codimensional in $X$ and denote by $c_i$ 
the codimension of $X_i$ in $X$. Then there is an isomorphism
of $A^*(pt)$-modules
$$
\bigoplus_{i=0}^{n} A^{k-c_i}(Y_i) \xra{\cong} A^k(X),
\quad\forall k\in \zz.   
$$
\end{thm}

The proof consists of the following steps. 
First, assuming that the localization sequences induced by the filtration
split, we obtain a decomposition isomorphism as stated in the theorem. 
Next, using the `elimination of points of indeterminacy'
by \cite{Hi} we produce a useful diagram (\ref{diagr}).
By means of this diagram we construct sections which split the
localization sequences. The properties of push-forwards with support
stated in Theorem \ref{relgysin} are crucial in the last step.

\begin{proof}{\bf I.}
For $0\leq i \leq n-1$, consider the localization sequence (cf. \ref{locseq})
$$
\ldots \ra A^*_{X_{i+1}}(X) \xra{s_i^*} A^*_{X_i}(X) \xra{u_i^*} 
A^*_{E_i}(X\smallsetminus X_{i+1}) \xra{\partial} A^{*+1}_{X_{i+1}}(X)
\ra\ldots,
$$
where $s_i^*$ is the extension of support and 
$u_i^*$ is induced
by the open embedding 
$u_i:(X\smallsetminus X_{i+1},X \smallsetminus X_i)
\ra (X,X\smallsetminus X_i)$.
Let $\tau_i$ denote the composition 
$$
A^{*}(Y_i)\xra{p_i^*} A^{*}(E_i) \xra{\imath_{th}} 
A^{*+c_i}_{E_i}(X\smallsetminus X_{i+1})
$$
which is an isomorphism by \ref{hominvar} and \ref{gysin}. 
Plugging these isomorphisms into the localization sequence 
we get
$$
\ldots \ra A^*_{X_{i+1}}(X) \xra{s_i^*} A^*_{X_i}(X) \xra{(\tau_i)^{-1}\circ u_i^*} 
A^{*-c_i}(Y_i) \xra{\partial\circ\tau_i} A^{*+1}_{X_{i+1}}(X)
\ra\ldots
$$

Assume that for each $i$ this exact sequence splits by means
of a section $\theta_i: 
A^{*-c_i}(Y_i) \ra A^*_{X_i}(X)$ of the morphism $(\tau_i)^{-1}\circ u_i^*$.
Then for each $i$ we have the decomposition
\begin{equation}\label{deciso}
(s_i^*,\theta_i): 
A^*_{X_{i+1}}(X) \oplus A^{*-c_i}(Y_i) \xra{\cong} A^*_{X_i}(X).
\end{equation}
Assembling these together,
we obtain an isomorphism
\begin{equation}\label{dcp}
\theta=(\theta_0',\ldots,\theta_n'): \bigoplus_{i=0}^n A^{*-c_i}(Y_i) 
\xra{\cong} A^*(X),
\end{equation}
where $\theta_i'=s_0^*\circ s_1^*\circ \ldots \circ s_{i-1}^*\circ \theta_i$.
It remains to construct sections $\theta_i$. 

\paragraph{II.}
Since any affine automorphism of $\af^r$
can be extended to an automorphism 
of $\pl^r=\pl(\af^r\oplus 1)$, we can consider the projective
completion $\pl(E_i\oplus 1)$
of the affine bundle $E_i/Y_i$ 
the same way we do it for vector bundles.

Observe that $E_i$ is an open subset of both $\pl(E_i\oplus 1)$ and $X_i$.
Thus $\pl(E_i\oplus 1)$ is birationally isomorphic to the closure
of $E_i$ in $X_i$ (which may not coincide with $X_i$ if the latter
is not irreducible.)

 
Then by results of Hironaka \cite[Cor.~3]{Hi} 
``on elimination of points of indeterminacy''
there exists a sequence of blow-ups of $\pl(E_i\oplus 1)$ with smooth
centers outside of $E_i$, resulting in a variety $W_i$,
together with a morphism 
$f_i:W_i \ra X_i$ which is identity on $E_i$.
Hence, we get a commutative diagram of the form
\begin{equation}\label{diagr}
\xymatrix{
X_i & E_i \ar[l]_{q_i}\\
W_i \ar[u]^{f_i}\ar[dr]_{\bar{p}_i}& E_i \ar[l]_{g_i}\ar[u]_{\id}\ar[d]^{p_i}\\
  & Y_i
}
\end{equation}
where $g_i$ is the open embedding $E_i \subset W_i$ and 
the projective morphism $\bar{p}_i$ is the composite of the projection maps
$W_i \ra \pl(E_i\oplus 1)\ra Y_i$. 
It can be checked (using Zariski Main Theorem) 
that the square part of the diagram is Cartesian.
Observe that by construction $W_i$ is a smooth projective variety over $k$.

\paragraph{III.}
We claim that the composition 
$$
A^*(Y_i)\xra{\bar{p}_i^*} A^*(W_i) \xra{\overline{(j_if_i)}_*} A^{*+c_i}_{X_i}(X),
$$
is the desired section $\theta_i$, where $j_i$ is the closed
embedding $X_i\hra X$. 
Note that we cannot consider the composition of push-forwards
$$
A^*(W_i)\xra{(f_i)_*} A^*(X_i) \xra{(j_i)_*} A^{*+c_i}(X)
$$
as $X_i$ is not necessarily smooth and $A^*(X_i)$ is not defined in general.
Instead, we fix a decomposition for $j_if_i$ as in \ref{relgysin} and consider
the associated push-forward with support $\overline{(j_if_i)}_*$.
By the first and third properties of push-forwards with support applied to
the square part of diagram (\ref{diagr}), 
we have a commutative diagram
(compare to (\ref{diagrpush}) with
$W=W_i$, $X=X$, $Z=X_i$, $U=X\smallsetminus X_{i+1}$, $W_U=E_i$, $Z_U=E_i$
and $f=j_if_i$) 
$$
\xymatrix{
A^{*+c_i}_{X_i}(X) \ar[r]^-{u_i^*} & A^{*+c_i}_{E_i}(X\smallsetminus X_{i+1}) \\
A^*(W_i) \ar[r]^{g_i^*} \ar[u]^{\overline{(j_if_i)}_*} & A^*(E_i) \ar[u]_{\imath_{th}}
}
$$
Combining this with the diagram
$$
\xymatrix{
A^*(W_i)\ar[r]^{g_i^*} & A^*(E_i) \\
  &  A^*(Y_i)\ar[ul]^{\bar{p}_i^*}\ar[u]_{p_i^*}
}
$$
we obtain
$$
u_i^* \circ\theta_i=(u_i^*\circ \overline{(j_if_i)}_* )\circ \bar{p}^*_i= \imath_{th}\circ (g_i^*\circ \bar{p}^*_i) = \imath_{th}\circ p_i^*=\tau_i.
$$
This completes the proof of Theorem \ref{decthm}.
\end{proof}

\begin{rem} The restriction on characteristic of the base field 
is due to the fact that the
resolution of singularities is required in the proofs of \cite{Hi}. 
\end{rem}

\section{Cellular decomposition in terms of \\$A$-correspondences}
\label{cellcorrdec}

We define the category of $A$-correspondences
following the strategies of \cite{Ma}.

\begin{ntt} The {\it category of $A$-correspondences}, denoted by $\corr$,
is the category whose
objects are smooth projective varieties over
$k$ and the set of morphisms from $Y$ to $X$
is defined by $\corr (Y,X) = A(Y\times X)$. 
An element of the ring $A(Y\times X)$ is called a {\it correspondence} between 
$Y$ and $X$.
For any $\alpha\in \corr(Y,X)$ and $\beta\in \corr(Z,Y)$ the correspondence
\begin{equation}\label{compcorr}
\alpha\circ \beta = (p_{13})_*(p_{12}^*(\beta)\cup p_{23}^*(\alpha)) \in \corr(Z,X)
\end{equation}
is the {\it composition} of $\alpha$ and $\beta$
(here $p_{12},\, p_{13},\, p_{23}$ denote the projections of
$Z\times Y\times X$ to $Z\times Y$, $Z\times X$, and $Y\times X$ respectively).
Observe that the push-forward $(p_{13})_*$ exists, since the cohomology theory $A$ is oriented
(see \ref{integr}). 
\end{ntt}

\begin{ntt}
The image of $\alpha\in\corr(Y,X)$ under the twisting map
$A(Y\times X)\ra A(X\times Y)$ is called the {\it transpose} of $\alpha$
and is denoted by $\alpha^t$. Clearly, $\alpha^t\in\corr(X,Y)$ and
$(\alpha^{t})^t=\alpha$. 
This operation makes $\corr$ a self-dual category.
\end{ntt}

\begin{ntt}\label{prodcor}
For any $\alpha\in \corr(Y,X)$ and $\beta\in \corr(Y',X')$ the correspondence
$$
\alpha\times\beta =p_{13}^*(\alpha)\cup p_{24}^*(\beta)\in \corr(Y\times Y',X\times X'),
$$
is called the {\it product} of $\alpha$ and $\beta$ (here
$p_{13}$ and $p_{24}$ denote the projections of $Y\times Y'\times X\times X'$
to $Y\times X$ and $Y'\times X'$ respectively).
Observe that the product of varieties and correspondences induces 
a tensor product structure on $\corr$. 
\end{ntt}

\begin{ntt}\label{graph}
Let $f: X\ra Y$ be a morphism of smooth projective varieties.
Its {\it graph} is the morphism
$$
\Gamma_f: X \xra{(f,\id)} Y\times X\, ,
$$
which is a closed embedding.
Let $c(f)=(\Gamma_f)_*(1_{A(X)}) \in A(Y\times X)$.
The assignment $X\mapsto c(X)=X$ and $f\mapsto c(f)\in\corr(Y,X)$
defines a contravariant functor $c: \smpr^{op}\ra\corr$ from the category
of smooth projective varieties to the category of correspondences.
\end{ntt}

\begin{ntt}
A crucial fact about the category $\corr$ is that the functor $A$ restricted 
to projective varieties factors through it. For  a correspondence
$\alpha\in\corr(Y,X)$,
define its realization $A(\alpha): A(Y)\ra A(X)$ as follows. 
Identify $A(Y)$ with 
$A(pt\times Y)=\corr(pt,Y)$ and assign to each $\beta\in A(Y)$ the composition
$\alpha\circ \beta\in\corr(pt,X)=A(X)$. This yields a covariant functor
$\corr\ra \ab$, which we will also denote by $A$. By (\ref{compcorr}),
$A(\alpha)(\beta)=(p_X)_*(p_Y^*\beta \cup \alpha)$.

Using the projection formula (see \ref{projform}), we check that 
$A(c(f))=f^*$ in the notation of \ref{graph}. 
Thus we have a diagram of functors
$$
\xymatrix{
\smpr^{op} \ar[rr]^A \ar[dr]_c& & \ab^{\zz} \\
    & \corr \ar[ur]_A& 
}
$$
\end{ntt}

\begin{ntt}\label{yoneda}
By the Yoneda Lemma, a correspondence $\alpha \in A(Y,X)$ 
is an isomorphism in $\corr$
if and only if the induced map 
$$
A(\id_T\times \alpha): A(T\times Y) \ra A(T\times X)
$$
is an isomorphism for any smooth projective variety $T$.
\end{ntt}

\begin{ntt}\label{grad}
Since $A$ is graded, we may endow morphisms in $\corr$ with grading
by setting 
$$\corr^c(Y,X)=A^{\dim Y +c}(Y\times X).$$
We say a correspondence $\alpha$ has degree $c$ if it is an element of
the group $\corr^c(Y,X)$. Observe that the transpose $\alpha^t$
has degree $\dim Y + c -\dim X$. 
\end{ntt}

In the sequel we will need the following fact.
\begin{lem}\label{coinc}
Let $h: W \ra Y\times X$ be a morphism
of smooth projective varieties. 
Put $\alpha=h_*(1_W)\in A(Y\times X)=\corr(Y,X)$. 
Define $pr_X$ and $pr_Y$ by the diagram 
$$
\xymatrix{
Y & Y\times X \ar[l]_{p_Y}\ar[r]^{p_X}& X \\
    & W \ar[u]^h\ar[ul]^{pr_Y}\ar[ur]_{pr_X}& 
}
$$
Then $(pr_X)_*\circ (pr_Y)^*=A(\alpha): A(Y) \ra A(X)$.
\end{lem}

\begin{proof}
We have
$$
(pr_X)_*\circ (pr_Y)^*=(p_X\circ h)_*\circ (p_Y\circ h)^*=
(p_X)_*\circ (h_* \circ h^*) \circ (p_Y)^*, 
$$
and projection formula \ref{projform} implies
$$
h_*\circ h^*(\gamma)=h_*(h^*(\gamma)\cup 1_W)=\gamma\cup h_*(1_W)=\gamma \cup \alpha, \quad \forall \gamma \in A(Y\times X).
$$
This proves the lemma.
\end{proof}

Now we are ready to state and prove the main result of this section.

\begin{thm}\label{corrdcp}
Let $A^*$ be an oriented cohomology theory over a field of
characteristic 0. Let $X$ be a relative cellular
space with a given filtration 
$X=X_0\supset X_1 \supset \ldots \supset X_n$. Denote by $r_i$ 
the rank of the affine bundle $E_i=X_i\smallsetminus X_{i+1}$ over $Y_i$. 
Then there is an isomorphism
of $A^*(pt)$-modules
$$
A^k(X)\xra{\cong} \bigoplus_{i=0}^{n} A^{k-r_i}(Y_i),
\quad\forall k\in \zz.   
$$
\end{thm}

The key idea of the proof is the following observation.
Since the section $\theta_i$ 
(constructed by means of diagram~(\ref{diagr}))
is the composite of a pull-back and a push-forward, 
it can be lifted to the category of $A$-correspondences.
In particular, the isomorphism~(\ref{dcp}) can be viewed
as the realization of some isomorphism $\alpha$ 
in the category of $A$-correspondences. 
As $\alpha$ is an isomorphism in $\corr$, 
its transpose $\alpha^t$ is an isomorphism in $\corr$ as well. 
Taking its realization $A(\alpha^t)$
we get precisely the desired isomorphism of cohomology groups of Theorem 
\ref{corrdcp}.

\begin{proof}
We want to show that the isomorphism $\theta$ of (\ref{dcp}) 
arises from some correspondence $\alpha$
via the functor $A: \corr \ra \ab^{\zz}$. Recall diagram (\ref{diagr}) and 
define 
the projective morphism 
$$
h_i=(\bar{p}_i,j_if_i): W_i \ra Y_i\times X,
$$ 
where $j_i:X_i\hra X$ is the closed embedding.
Define $\alpha_i:=(h_i)_*(1_{W_i})$.
As $\codim h_i=\dim (Y_i\times X) -\dim W_i = \dim Y_i + c_i$,
we conclude that 
$$
\alpha_i\in A^{\dim Y_i +c_i}(Y_i\times X)=\corr^{c_i}(Y_i,X).
$$
We have 
$$
A(\alpha_i)=(j_if_i)_*\circ\bar{p}_i^*=(s\circ \overline{(j_if_i)}_*)\circ 
\bar{p}_i^*=s\circ\theta_i=\theta_i',
$$
where the map $s=s_0^*\circ\ldots\circ s_{i-1}^*$ is the extension of support $A_{X_i}(X)\ra A(X)$, 
the first equation holds by Lemma \ref{coinc} and the second follows
from the second property of push-forwards with support (see \ref{relgysin}).
Now put $\alpha =(\alpha_0,\ldots,\alpha_n)\in \corr(\amalg Y_i,X)=\corr(Y,X)$.
Then $A(\alpha)=\theta$.

\begin{prop}
The correspondence $\alpha$ is an isomorphism in $\corr$.
\end{prop}

\begin{proof}
By Yoneda's Lemma \ref{yoneda} it suffices to check that for any $T\in\corr$ the map
$\alpha_T: \corr(T,Y)\ra \corr(T,X)$, given by $\alpha_T(\beta)=\alpha\circ\beta$
for $\beta\in \corr(T,Y)$, is bijective. A straightforward computation shows that
$\alpha_T=A(1_T\times\alpha)$, where the product of correspondences $1_T$ and $\alpha$
is defined in \ref{prodcor}. Thus we must check that 
$A(1_T\times \alpha): A(T\times Y)\ra A(T\times X)$ is an isomorphism 
for any smooth projective variety $T$.
Another straightforward verification shows that we can insert $(T\times-)$ throughout
the proof of Theorem \ref{decthm} and 
\begin{enumerate}
\item[(i)] get a cellular decomposition for $T\times X$ with bases $T\times Y_i$;
\item[(ii)] get splittings $\theta_i$ for $T\times X$ the same way it was done
for $X$ by multiplying diagram (\ref{diagr}) with $T$;
\item[(iii)] get an isomorphism $\bigoplus_{i=0}^n A^{*-c_i}(T\times Y_i) \ra A^*(T\times X)$;
\item[(iv)] check that this isomorphism arises from the correspondence $1_T\times\alpha$.
\end{enumerate}
This proves the proposition.
\end{proof}

Consider now $\alpha^t=(\alpha_0^t,\ldots,\alpha_n^t)\in\corr(X,\amalg Y_i)=\corr(X,Y)$.
Note that $\deg \alpha_i^t= \dim Y_i+ c_i -\dim X =-r_i$ by \ref{grad}.
Being the transpose of an isomorphism, $\alpha^t$ must be an isomorphism as well.
Let $\sigma_i=A(\alpha_i^t)$. Then
$\sigma=(\sigma_0,\ldots,\sigma_n)=A(\alpha^t)$ 
is the desired isomorphism of Theorem \ref{corrdcp}.
The theorem is proved.
\end{proof}

\section{Decomposition of $A$-motives}\label{amot}

We define the category of $A$-motives following \cite{Ma}.

\begin{ntt}
Recall that an additive category $C$ is called pseudo-abelian
if for any projector $p\in Hom(X,X)$, $X\in C$, there exists a kernel
$\ker p$, and the canonical homomorphism $\ker p \oplus \ker(\id_X-p) \ra X$
is an isomorphism. Clearly, $\id_X -p$ is also a projector in $C$.

Let $C$ be an additive category. Its pseudo-abelian completion
is the category $\tilde{C}$ defined as follows.
The objects are pairs $(X,p)$, where $X$ is an object of $C$ and
$p\in Hom(X,X)$ is a projector. The morphisms are given by
$$
Hom_{\tilde{C}}((X,p),(Y,q)) = \frac{\{\alpha\in Hom(X,Y)\mid \alpha\circ p=q\circ \alpha\}}{\{\alpha \in Hom(X,Y)\mid \alpha\circ p = q\circ \alpha=0\}}.
$$ 

If $p_0, p_1, \ldots, p_n \in Hom(X,X)$ are projectors satisfying
$p_i\circ p_j=0$, $\forall i\neq j$, and 
$\sum_{i=0}^n p_i=\id_X$,
then in $\tilde{C}$ we have
$$
(X,\id_X)\cong \bigoplus_{i=0}^n (X,p_i).
$$
To simplify the notation we will write $X$ for $(X,\id_X)$.
\end{ntt}

\begin{ntt}
Consider the category $\corr^0$ of correspondences
of degree $0$.
The pseudo-abelian completion of this category 
will be called the {\it category
of $A$-motives} and will be denoted by $\mot$.
The following properties of this category can be proved by a 
standard reasoning (or see \cite{Ma}).

(i) There is a covariant functor $\corr^0 \ra \mot$, $X\mapsto X=(X,\id_X)$.

(ii) The cohomology theory $A^*$ factors through $\mot$.

(iii) There is a tensor product structure on $\mot$ induced by the
tensor product structure on $\corr$.
\end{ntt}

Projective Bundle Theorem \ref{projbun}
for $A$ implies the following

\begin{prop}
Let $E\ra X$ be a vector bundle of rank $n+1$ over a smooth projective
variety $X$. Then there exist $n+1$ orthogonal projectors
$p_0,p_1,\ldots,p_n\in\corr^0(\pl(E),\pl(E))$ 
(defined by particular formulas and playing non-interchangeable roles)
such that $\sum_{i=0}^n p_i = \id_{\pl(E)}$, i.e., in $\mot$ one has
a decomposition
$$
\pl(E)\cong \bigoplus_{i=0}^n (\pl(E), p_i).
$$
The terms of this decomposition can be written as
$$
(\pl(E),p_i)\cong X\otimes L^{\otimes i},\quad 0\leq i\leq n
$$ 
where $L=(\pl^1,p_1)$ is called the Lefschetz motive.
Observe that $L$ exists in the category $\mot$ but not in $\corr$.
\end{prop}

\begin{proof}
See \cite[Sect.~7]{Ma}.
\end{proof}

\begin{ntt}\label{motgrad}
For the groups of morphisms in $\mot$ we have the isomorphisms
$$
Hom_{\mot}(X,Y)\cong Hom_{\mot}(X\otimes L^i, Y\otimes L^i), 
\quad \forall i\geq 0.
$$
Hence, one can endow morphisms in $\mot$ with grading that respects
the grading in $\corr$ by setting
$$
Hom_{\mot}^c(X,Y)=\varinjlim_i Hom_{\mot}(X\otimes L^{c+i},Y\otimes L^i), \quad \forall c\in \zz.
$$
It can be checked that a correspondence $\alpha\in Hom_{\corr}^c(Y,X)$ 
of degree $c\geq 0$
induces a morphism 
$$\alpha\in Hom_{\mot}^c(Y,X)=Hom_{\mot}(Y\otimes L^c,X)$$
of degree $c$ in the category of motives, see \cite[Sect.~8]{Ma}.
If $\alpha^t\in Hom_{\corr}^{-r}(X,Y)$, where $r=\dim X -\dim Y -c$, is the transpose
of $\alpha$, then it induces the morphism 
$$\alpha^t\in Hom_{\mot}^{-r}(X,Y)=Hom_{\mot}(X,Y\otimes L^r).$$ 
\end{ntt}

Now we are ready to state and prove the main result of this paper.

\begin{thm}\label{mainthm}
Let $A^*$ be an oriented cohomology theory over a field
of characteristic 0 and 
$X$ a relative cellular space with a filtration
$X=X_0\supset X_1 \supset\ldots\supset X_n\supset\emptyset$.
Then there exists an isomorphism in the category of $A$-motives
$$
X\xra{\cong} \bigoplus_{i=0}^n Y_i\otimes L^{r_i},
$$
where $L$ is the Lefschetz motive and $Y_i$ is the motive of the base 
of the affine bundle
$p_i:X_i\smallsetminus X_{i+1}\ra Y_i$ of rank $r_i$.
\end{thm}

\begin{proof}
Consider the correspondence $\alpha=(\alpha_0,\ldots,\alpha_n)$ defined
in the proof of \ref{corrdcp}. Recall that $\alpha_i\in\corr^{c_i}(Y_i,X)$.
According to \ref{motgrad} it induces the morphism in the category $\mot$ 
$$
\alpha=(\alpha_0,\ldots,\alpha_n)\in
Hom_{\mot}(Y, X),
$$
where $Y=\bigoplus_i Y_i\otimes L^{c_i}$ in $\mot$, and
its realization $A(\alpha)$ coincides with the isomorphism 
of Theorem~\ref{decthm}.

Since $\alpha$ is an isomorphism in $\corr$ it is an isomorphism
in $\mot$ as well.
Take now its transpose $\alpha^t \in Hom_{\mot}(X, Y^t)$,
where $Y^t=\bigoplus_i Y_i\otimes L^{r_i}$ (see \ref{motgrad} about
degrees conventions). Clearly, 
it gives the desired isomorphism in $\mot$
$$
\alpha^t: X \xra{\cong} \bigoplus_{i=0}^n Y_i\otimes L^{r_i}
$$
In particular, its realization $A(\alpha^t)$ gives the
decomposition of cohomology groups of Theorem \ref{corrdcp}.
The theorem is proved.
\end{proof}

\section{Relative cellular spaces with smooth filtrations}

In the present section we provide examples of relative
cellular spaces with {\it smooth
filtrations} (i.e., filtrations with smooth terms), as opposed
to most examples of \cite{CGM} in which the filtrations are not smooth.
Observe that it is easier to deduce formulas (\ref{cdecomp}) and 
(\ref{rdecomp}) for smooth filtrations than in general, for
one does not need push-forwards with support in this case.
Working with smooth filtrations, 
one may expect to apply the arguments of the present paper
to derived Witt groups as soon as one has push-forwards for them. 
Some results in this direction
have been obtained recently in \cite{CH} and \cite{Ne}.

We consider the cases of completely split quadrics \ref{quadr}
and Grassmannians \ref{grass}. In these cases it is possible to
construct the varieties $W_i$ of (\ref{diagr}) explicitly, without
resolution of singularities. This implies 
that Corollaries~\ref{motquad} and \ref{motgras} hold 
over any field of characteristic different
from 2.

\begin{ntt}[Completely split quadrics]\label{quadr}
Let $Q=(d+1)\hh$ be a split projective quadric of dimension $2d$.
Then $Q$ can be thought of as the hypersurface in $\pl^{2d+1}$ given by
$$
Q=\{[x_0:y_0:\ldots :x_d:y_d]\in \pl^{2d+1}\mid \sum_{i=0,d}x_iy_i=0\}.
$$
We claim that $Q$ has a smooth filtration defined as follows.

Consider two linear subspaces of $\pl^{2d+1}$ of dimension $d$ defined by the equations
$$
Z_1=\{[x_0:y_0:\ldots :x_d:y_d]\in \pl^{2d+1}\mid x_0=x_1=\ldots =x_d=0\}\quad\text{and}
$$
$$
Z_2=\{[x_0:y_0:\ldots :x_d:y_d]\in \pl^{2d+1}\mid y_0=y_1=\ldots =y_d=0\}.
$$
Clearly $Z_1 \cong Z_2 \cong \pl^d $ and $Z_1\cap Z_2=\emptyset$. Consider the linear map
$\pi_1:\pl^{2d+1}\smallsetminus Z_2 \ra Z_1$ given by the projection to the $y$-coordinates.
Both subspaces $Z_1$ and $Z_2$ lie in $Q$. It can be easily seen that
the restriction $\pi_1: Q\smallsetminus Z_2\ra Z_1$  
naturally has a structure of a vector bundle of rank $d$.
Hence, we get a smooth filtration on $Q$: 
$$X_0=Q \supset X_1=Z_2, \qquad Y_0=Z_1, Y_1=Z_2.$$
\end{ntt}

Theorem \ref{mainthm} applied to this filtration yields the following

\begin{thm}\label{motquad}
Let $X=(d+1)\hh$, $d\geq 0$, be a split projective quadric of dimension $2d$.
Then in the category of $A$-motives
$$
X \cong \pl^d \oplus (\pl^d\otimes L^d).
$$
In particular,
$$
A^k((d+1)\hh)\cong A^k(\pl^d)\oplus A^{k-d}(\pl^d),\quad \forall k\in \zz
$$
For the Chow groups $A^*=CH^*$ it gives the well-known decomposition
$$
CH^*((d+1)\hh)\cong \bigoplus_{i=0}^d (\zz [i] \oplus \zz [i+d]),
$$
where $\zz[i]$ is the abelian group $\zz$ put in degree $i$.

For the $K$-functor, i.e., $A^k(X,U)=K_{-k}(X,U)$ 
(see \cite[Ex.~2.1.8]{Pa})
it gives the decomposition
$$
K_0((d+1)\hh)\cong \zz^{2d+2}
$$
\end{thm}

\begin{ntt}[Grassmannians]\label{grass}
Let $Gr(d,n)$ be the Grassmann variety of $d$-planes in the $n$-dimensional
affine space $\af^n$. Then $Gr(d,n)$ has a smooth
filtration defined as follows.

Fix a $(n-1)$-dimensional linear subspace $V$ of $\af^n$.
Then there is the induced closed embedding 
$Gr(d,n-1)\ra Gr(d,n)$ of smooth projective
varieties, where $Gr(d,n-1)$ is identified with the
variety of $d$-planes in $V$. 
Consider the complement $E=Gr(d,n)\smallsetminus Gr(d,n-1)$.
We claim that this is a vector bundle over $Gr(d-1,n-1)$.
Indeed, an element of $E$ is a $d$-plane that intersects $V$
along a linear subspace of dimension $d-1$, i.e., it gives
a point of $Gr(d-1,n-1)$.
Moreover, the set of all elements of $E$ that intersect $V$
along a fixed linear subspace of dimension $d-1$ can be identified
with the affine space of dimension $n-d$.
This gives a vector bundle structure.

Applying this procedure recursively
we get a smooth cellular filtration
$$
X_0=Gr(d,n)\supset X_1=Gr(d,n-1) \supset \ldots \supset X_{n-d}=
Gr(d,d)= pt
$$
with the bases $Y_0=Gr(d-1,n-1),Y_1=Gr(d-1,n-2),\ldots, Y_{n-d}=pt$. 
\end{ntt}

Theorem \ref{mainthm} recursively applied to this filtration yields the 
following assertion which is well-known for Chow groups \cite{Fu}.

\begin{thm}\label{motgras}
Let $X=Gr(d,n)$ be the Grassmann variety of $d$-planes in the affine space
$\af^n$ of dimension $n$. Then there is the decomposition of $A$-motives
$$
X\cong \bigoplus_{\lambda} L^{|\lambda|},
$$
where the sum is taken over all partitions 
$\lambda=(\lambda_1,\ldots,\lambda_d)$ with 
$1\leq \lambda_1< \ldots <\lambda_d\leq n$ and $|\lambda|=\lambda_1+\ldots+\lambda_d$.

In particular,
$$
A^k(Gr(d,n))\cong \zz^{\#\{\lambda\mid k=|\lambda|\}}.
$$
\end{thm}

\noindent
Alexander Nenashev\\
Mathematics Department\\
Glendon College (York University)\\
2275 Bayview Av.\\
Toronto, ON\\
Canada M4N 3M6\\
e-mail: \href{mailto:Nenashev@gl.yorku.ca}{Nenashev@gl.yorku.ca}

\ 

\noindent
Kirill Zainoulline\\
Fakult\"at f\"ur Mathematik\\
Universit\"at Bielefeld\\
Postfach 100131\\
D-33501 Bielefeld\\
Germany\\
e-mail: \href{mailto:kirill@mathematik.uni-bielefeld.de}{kirill@mathematik.uni-bielefeld.de}


\begin{thebibliography}{999}


\bibitem[Br]{Br} Brosnan,~P.\;
On motivic decomposition arising from the method of Bialynicki-Birula.
Preprint Server 
``{\it Linear Algebraic Groups and Related Structures}'' (2004), 
no.147, 20 pp.\\
\url{http://www.mathematik.uni-bielefeld.de/LAG/man/147.html}

\bibitem[CGM]{CGM} Chernousov, V.; Gille, S.; Merkurjev, A.\;
Motivic decomposition of isotropic projective homogeneous varieties.
(to appear in {\it Duke Math. J.}) Preprint (2003), 17 pp.\\
\url{http://www.math.ucla.edu/~merkurev/index.html}

\bibitem[CH]{CH} Calm\`es, B.; Hornbostel, J.\;
Witt Motives, Transfers and Reductive Groups. 
Preprint Sever 
``{\it Linear Algebraic Groups and Related Structures}'' (2004),
no.143, 29 pp.\\
\url{http://www.mathematik.uni-bielefeld.de/LAG/man/143.html}

\bibitem[Fu]{Fu} Fulton, W.\;
Intersection theory. Second edition. 
Ergebnisse der Mathematik und ihrer Grenzgebiete. 3. Folge. 
A Series of Modern Surveys in Mathematics, 2. {\it Springer-Verlag}, 
Berlin, 1998. xiv+470 pp. 

\bibitem[Gr]{Gr} Grothendieck, A.\;  
Sur quelques propri\'et\'es fondamentales en th\'eorie 
des intersections. S\'eminaire C.~Chevalley: 
anneaux de Chow et applications. {\it Secr\'etariat
math\'ematique}, Paris 1958.

\bibitem[Hi]{Hi} Hironaka, H.\;
Flattening theorem in complex algebraic geometry. 
{\it Amer. J. of Math.} 97 (1975), 503--547.

\bibitem[LM]{LM} Levine, M.; Morel, F.\;  
Cobordisme alg\'ebrique I. 
{\it C. R. Acad. Sci. Paris, S\'er. I Math.} 332 (2001), no.8, 723--728. 

\bibitem[Ka]{Ka} Karpenko, N.\; 
Cohomology of relative cellular spaces and of isotropic flag varieties. 
{\it St.-Petersburg Math. J.} 12 (2001), no.1, 1--50.

\bibitem[Ma]{Ma} Manin, Ju.\;
Correspondences, Motifs and Monoidal Transformations.
{\it Math. USSR Sbornik} 6 (1968), 439--470.


\bibitem[Ne]{Ne} Nenashev, A.\;
On the Witt groups of projective bundles
and split quadrics: Geometric reasoning. 
{\it K-theory Preprint Archives} (2004), no.696, 9 pp.
\url{http://www.math.uiuc.edu/K-theory/0696/}

\bibitem[Ne1]{Ne1}Nenashev, A.\; 
Gysin maps in oriented theories. 
{\it K-theory Preprint Archives} (2004), no.695, 9 pp.\\
\url{http://www.math.uiuc.edu/K-theory/0695/}

\bibitem[Pa]{Pa} Panin,~I. (after Panin, I. and Smirnov, A.)\; 
Oriented Cohomology Theories of Algebraic Varieties.
{\it K-Theory J.} 30 (2003), 265--314.

\bibitem[Pa1]{Pa1} Panin,~I. (after Panin, I. and Smirnov, A.)\; 
Push-forwards in oriented cohomology theories of algebraic varieties II.
{\it K-theory Preprint Archives} (2003), no.619, 83 pp.\\
\url{http://www.math.uiuc.edu/K-theory/0619/}

\bibitem[Sm]{Sm} Smirnov, A.\;
Oriented structures on cohomology theories
of algebraic varieties. Preprint (2004).

\end{thebibliography}
\end{document}